\documentclass[a4paper,11pt]{article}
\usepackage{amssymb,delarray,amsmath}

\def\cc{{\mathbb C}}

\def\rr{{\mathbb R}}
\def\nn{{\mathbb N}}

\def\sph{{\mathbb S}}

\def\st{_{\mathsf{st}}}
\def\eps{\varepsilon}

\def\st{_{\mathsf{st}}}

\def\lim{\mathop{\sf{lim}}}

\def\eps{\varepsilon}

\def\<{\langle}\let\la=\<
\def\>{\rangle}\let\ra=\>

\def\d{\partial}

\def\ddef{\mathrel{{=}\raise0.3pt\hbox{:}}}
\def\deff{\mathrel{\raise0.3pt\hbox{\rm:}{=}}}

\def\fraction#1/#2{\mathchoice{{\msmall{ #1\over#2}}}%
{{ #1\over #2 }}{{#1/#2}}{{#1/#2}}}

\def\le{\leqslant}

\def\longpoints{\leaders\hbox to 0.5em{\hss.\hss}\hfill \hskip0pt}
\def\stateskip{\smallskip}
\def\state#1. {\stateskip\noindent{\bf#1. }} %\medskip
\def\statep#1. {\stateskip\noindent{\bf#1 }} %\medskip

\def\proof{\state Proof. \2}

\def\Chi{\raise 2pt\hbox{$\chi$}}

\def\ie{\hskip1pt plus1pt{\sl i.e.\/,\ \hskip1pt plus1pt}}

\def\sli{{\sl i)} } 
\def\slii{{\sl i$\!$i)} } 
\def\sliii{{\sl i$\!$i$\!$i)} }

\def\Chi{\raise 2pt\hbox{$\chi$}}

\let\phI=\phi\let\phi=\varphi\let\varphi=\phI
\def\eps{\varepsilon}

\def\d{\partial}

\def\1{{1\mkern-5mu{\rom l}}}

\def\ge{\geqslant}

\def\fraction#1/#2{\mathchoice{{\msmall{ #1\over#2}}}%
{{ #1\over #2 }}{{#1/#2}}{{#1/#2}}}

\def\le{\leqslant}

\newcommand{\2}{\thinspace}

\def\qed{\ \ \hfill\hbox to .1pt{}\hfill\hbox to .1pt{}\hfill $\square$\par}

\def\comment#1\endcomment{}

\def\eqqno(#1){\label{(#1)}}
\def\eqqref(#1){(\ref{(#1)})}

\title{On nonimbeddability of Hartogs figures into complex manifolds}
\author{E. Chirka \and S. Ivashkovich}
\begin{document}
\maketitle
\newtheorem{lem}{Lemma}
\newtheorem{quest}{Question}
\newtheorem{cor}{Corollary}

\paragraph{\S 1. Introduction.} In this note we whant to discuss
two closely related questions about possibility of certain
imbeddings  of Hartogs figures into general complex manifolds. In
fact, we shall construct a counterexample to both these questions.

First question was asked by Evgeny Poletsky. Let $\Delta$ denote the 
unit disk in $\cc$, $\Delta(r)$ - disk of radius $r$, $\Delta^2$ - unit 
bidisk in $\cc^2$ and $\sph^1$ - the unit cirle, $A^{r_2}_{r_1}$ - an annulus 
$\Delta_{r_2}\setminus \bar\Delta_{r_1}, r_1<r_2$. Recall that:

\smallskip\noindent (a) a "thin Hartogs figure" is the following
set in $\cc^2$
\begin{displaymath}
H = \{ (z,w)\in \cc^2: z=0, |w|\le 1 \textrm{ or } |\textrm{Re}\,z|\le 1, 
\textrm{Im}\,z=0, |w|=1\} = 
\end{displaymath}
\begin{displaymath}
= (\{0\}\times\bar\Delta)\cup ([-1,1]\times\sph^1);
\end{displaymath}

\smallskip
\noindent (b) a "thick Hartogs figure" (or simply Hartogs figure) is a set
of the form
\begin{displaymath}
H_{\eps}:=\{ (z,w)\in \cc^2: |z|<\eps , |w|<1+\eps
 \textrm{ or }  |z|<1+\eps, 1-\eps <|w|<1+\eps \} =
\end{displaymath}
\begin{displaymath}
=\Delta_{\eps}\times \Delta_{1+\eps}\cup \Delta_{1+\eps}\times A_{1-\eps}^
{1+\eps}, 
\end{displaymath}
\noindent for some $\eps$, $0<\eps <1$.

\smallskip 
Let $X$ be some complex manifold and suppose we are given a
continuous map $f:H\to X$ such that $f(0,w)$ is holomorphic 
on the disk $\{0\}\times \Delta $.

\begin{quest}{\bf.} Assume in addition that $f:H\to X$ is an 
imbedding. Can one find a neighborhood $V\supset f(H)$ which is 
biholomorphic to an open set in $\cc^n$ or, more generally, in 
some Stein manifold?
\end{quest}

\smallskip To formulate the second question suppose now that $X$
is of complex dimension two and is foliated by complex curves over
the unit disk. More precisely, 
there is a holomorphic submersion $\pi :X\to \Delta $ with connected 
fibers $X_z:=\pi^{-1}(z)$. Hartogs figures $H_{\eps }$ are naturally 
foliated
over the disk in the first factor $\cc_z$ of $\cc^2_{z,w}$ and we 
denote the corresponding projection by $\pi_1:H_{\eps}\to \Delta$.
A holomorphic mapping $f:(H_{\eps},\pi_1)\to (X,\pi )$ is called 
foliated if there exist a holomorphic map $\zeta :\Delta\to \Delta $
such that $\pi(f(z,w))=\zeta(z)$.

Suppose further we are given a
smooth family $\Gamma =\{ \gamma_z:z\in \Delta \}$ of diffeomorphic 
images of a circle $\sph^1$ with $\gamma_z\subset X_{z}$ such that 
$\gamma_0$ doesn't bound a disk in $X_0$
but there are $z\in \Delta $ arbitralily close to $0$ such that
$\gamma_{z}$ bound a disk in $X_z$. Denote by $\sph^1_a:=\{a\}\times 
\{|w|=1\}\subset \cc^2$ circles in corresponding fibers of $H_{\eps}$.  

\begin{quest}{\bf.} Does there exist $\eps >0$ 
such that a "thick" Hartogs figure $H_{\eps}$ can be holomorphically 
imbedded into $X$ in the following way:

\smallskip
\noindent 1) imbedding $f:H_{\eps}\to X$ is foliated;

\smallskip
\noindent 2) $f(\{0\}\times \Delta )\subset X_{a}$, $f(\sph^1_{0})$
homologous to $\gamma_{a}$  for some $a\in \Delta$;

\smallskip\noindent 3)  the curve $f(\sph^1_1)$ is contained in $X_0$ and 
is homologous to $\gamma_0$ in $X_0$?
\end{quest}

\smallskip In [Br-1] and [Br-2] the existence of such imbedding is 
used as an obvious fact, see p.124 and p.146 correspondingly. 

\smallskip The goal of this note is to provide an example giving the
negative answer to both Questions 1 and 2. We shall construct the
following

\medskip\noindent\bf
Example. {\it There exists a complex surface $X$  with a
holomorphic submersion $\pi$ onto the unit disk $\Delta $ such
that:

\smallskip
1) all fibers $X_z:=\pi^{-1}(z)$ are disks with possible punctures;

2) the fiber $X_0$ over the origin is a punctured disk; the subset
$U \subset \Delta$ consisting of such $z$ that the fiber $X_z$ is a
disk, is nonempty, open and $\partial U\ni 0$;

3) for any circle $\gamma_0$ around the puncture in $X_0$ and for
any circle $\gamma_{a}$ in any of $X_{a}, a\in U$, there does not exist
a foliated holomoprhic map $f$ from any "thick" Hartogs figure
$H_{\varepsilon}$ to $X$ such that $f(\{0\}\times\Delta)\subset
X_{a}$ and $f(\sph^1_1)\subset X_0$ is homologous to $\gamma_0$ in the 
fiber $X_0$.}

\smallskip\rm On the way of constructing this counterexample to the Question 2
we construct also a counterexample to the Question 1, which occurs to be
somewhat simpler (but not essentially simpler). We whant to emphasize that 
while existence of a counterexample to Q.2 should be of no surprise, the
existence of it to Q.1 is somewhat unexpectable, because the "thin Hartogs
figure" has very fiew of a complex structure - just one complex disk.

\smallskip\noindent\bf Asknowledgments. \rm We would like to thank 
E. Poletsky,
who was the first who asked us the question about possibility of imbeddings 
of Hartogs figures into general complex manifolds. We would like also 
to asknowledge M. Brunella for sending us a preprint [Br-3] where his 
erroneous argument with Hartogs figures is replaced by another approach
using a sort of "nonparametrized" Levi-type extension theorem.

\smallskip At any rate the question about possibility of certain imbeddings 
of Hartogs figures into a general complex manifold seems to be of growing
demand and interest.

\paragraph{\S 2. Construction of the example.}

\rm Our example is based on the violation of the argument principle.
Let $J\st$ denote the usual complex structure in $\cc^2_{z,w}$.

Take a function $\lambda(t)\in C^\infty(\rr), 0\le\lambda\le1,$
which satisfies

\begin{displaymath}
\lambda(t) = \left\{ \begin{array}{ll} 0 & \textrm{for $t<1/9$;}\\
1 & \textrm{for $t>4/9$}.
\end{array} \right.
\end{displaymath}

\noindent For $k\in \nn $ consider the following domain $M=M_k$ in
$\cc\times \Delta\subset \cc^2_{z,w}$: 
\begin{displaymath}
M:=(\cc\times\Delta)\setminus \{ (z,w): 
\frac{1}{3}\le |z|\le \frac{2}{3}, w^2=z^k\lambda(|z|^2) \textrm{ or }
|z|\ge \frac{1}{3}, w=0 \}.
\end{displaymath}

Let $J=J_k$ be the (almost) complex structure on $M_k$ with the basis of
(1,0)-forms costituted by $dz$ and $dw+a_kd\bar z$, where

\begin{displaymath}
a_k(z,w) = \left\{ \begin{array}{ll}
\frac{wz^{k+1}\lambda'(|z|^2)}
{w^2-z^k\lambda(|z|^2)} & \textrm{for $\frac{1}{3}<|z|<\frac{2}{3}$,}\\
0 & \textrm{otherwise.}
\end{array} \right.
\end{displaymath}

The subspace in $\Lambda^{p+q}(M)$ consisting of $(p,q)$-forms relative to 
$J$ we shall denote by $\Lambda_{J}^{p,q}(M)$.

\begin{lem}. $J$ is well defined on the whole of $M$, is
(formally) integrable, hence $(M, J)$ is a complex manifold.
Moreover

\sli $J=J\st$ on $M\setminus (\overline{A^{2/3}_{1/3}}\times 
\Delta)$;

\slii functions $f_k(z,w)=w+{z^k\over w}\lambda(|z|^2)$ and
$g(z,w)=z$ are $J$-holomorphic on $M$;

\sliii ${\sf ind}_{|w|=1-\eps}f_k(z,w)=-1$ for $|z|\ge 1$ and $0<\eps <1/6$.

\end{lem}
\proof \sli Integrability condition on $J$ reads as $d\Lambda_{J}^{1,0}
\subset \Lambda_{J}^{2,0} + \Lambda_{J}^{1,1}$ where $\Lambda_{J}^{2,0}$
is the linear span of $\Lambda_{J}^{1,0}\wedge \Lambda_{J}^{1,0}$ 
and $\Lambda_{J}^{1,1}$ is the same for $\Lambda_{J}^{1,0}\wedge 
\Lambda_{J}^{0,1}$. Any form $\alpha\in\Lambda_{J}^{1,0}$ is represented
as $\alpha_1 dz + \alpha_2(dw + a_kd\bar z)$ with smooth $\alpha_1, 
\alpha_2$, hence,  $d\alpha \equiv \alpha_2da_k\wedge d\bar z{\sf mod}
(\Lambda_j^{2,0}+\Lambda_j^{1,1})$.

Now, $da_k\wedge d\bar z=\frac{\d a_k}{\d z}dz\wedge d\bar z + 
\frac{\d a_k}{\d w}dw\wedge d\bar z + \frac{\d a_k}{\d \bar w}d\bar w\wedge
d\bar z\equiv \frac{\d a_k}{\d w}dw\wedge d\bar z$

\noindent ${\sf mod}\Lambda^{1,1}_J$
since $\frac{\d a_k}{\d \bar w}=0$. Finally, $dw\wedge \d \bar z=(dw + a_kd
\bar z)\wedge d\bar z\in \Lambda^{1,1}_J$.

\slii Really, $df_k(z,w)=(k\frac{z^{k-1}}{w}\lambda + \lambda'\frac{z^{k}}{w}
\bar z)dz + \lambda'\frac{z^{k+1}}{w}d\bar z + (1-\frac{z^k}{w^2}\lambda )dw=$

\noindent $=(k\frac{z^{k-1}}{w}\lambda + \lambda'\frac{z^{k}}{w}\bar z)dz + 
(1-\frac{z^k}{w^2}\lambda )(dw + a_kd\bar z)\in \Lambda^{1,0}_J$, \ie 
$\bar\d f_k=0$. The case of $g(z,w)=z$ is obvious since $dz\in \Lambda^{1,0}_J
(M)$.

\sliii is obvious.

\qed

\begin{lem}.
{\it Let $(\zeta,\eta)\mapsto(z(\zeta),w(\zeta,\eta))$ be any foliated
holomorphic map $(H_{\eps},J\st)\to (M,J)$ such that 

\sli $|z(0)|<1/3$, $w(0,0)=0$;

\slii $|w(\zeta,\eta)|\ge \delta $ for some $\delta>0$ and for all
$\zeta\in \Delta_{1+\eps} , |\eta |=1$.

\noindent Then $z(\Delta)\subset \Delta$.
}
\end{lem}
\proof \rm Suppose not. Set $U=z^{-1}(\Delta)$. Then $U\not= \Delta$
there are $\zeta_0\in \Delta\cap \d U$ and a curve $\gamma (t)$ 
from $\gamma (0)=z(0)$ to $\gamma (1)=\zeta_0$ such that $\gamma (t)\in U$ for 
$0\le t<1$. The function  $F_k(\zeta,\eta)=f_k(z(\zeta),w(\zeta,\eta))$
is holomorphic in $H_{\eps }$ and therefore holomorphically extends 
onto the bidisk $\Delta^2$.

Since $F_k(\zeta,\eta)=w(\zeta,\eta)+\frac{z(\zeta)^k}{w(\zeta,\eta)}
\lambda(|z(\eta)|^2)$, we see that ${\sf ind}_{|\eta |=1}F_k(0,\eta)=
{\sf ind}_{|\eta |=1}w(0,\eta)>0$ due to $J\st$-holomorphicity of 
$w(\zeta,\eta )$ on $\{ 0\}\times \Delta$. 

But $|z(\zeta_0)|=1$, so $\lambda (|z(\zeta_0)|^2)=1$ and therefore 
$\left| \frac{z(\zeta_0)^k}{w(\zeta,\eta)}\right| 
\lambda(|z(\zeta_0)|^2)>1>|w(\zeta_0,\eta)|$ for $|\eta|=1$. As $w(\zeta, \eta )$
is holomorphic on $\{ \zeta_0\}\times \Delta$, one has 
${\sf ind}_{|\eta |=1}F_k(\zeta_0,\eta)=
{\sf ind}_{|\eta |=1}\frac{1}{w(0,\eta)}<0$, which contradicts to the 
holomorphicity of $F_k$ on $\Delta^2$.

\qed

\smallskip\noindent
{\sl Counterexample to the Question 1.}
Define $f:H\to (M,J)$ as follows: $z(0,\eta )=0,
w(0,\eta)=(1-\eps )\eta$ on $\{0\}\times \bar\Delta$ and $z(\zeta,e^{i\theta})=
\zeta, w(\zeta,e^{i\theta})=
(1-\eps)e^{i\theta}$ with $0<\eps<1/10$ on $[-1,1]\times \sph^1$, \ie  $f$ is a 
scaled tautological imbedding.

\smallskip
\begin{lem} There is no neighborhood of $f(H)$ biholomorphic to an open set
in some Stein manifold.
\end{lem}
\proof Suppose that such a neighborhood $V\supset f(H)$ exists and
 $p:V\to Y$ is a biholomorphic imbedding of $V$ into a Stein manifold. 
 Let $\pi $ be the projection $(z,w)\to z$ of $M$ to $\cc$. After
shrinking $V$ if necessary we can assume about the projection $\pi|_V:V\to
\Delta $ the following:

\smallskip
\sli for $z$ in a neighborhood $W_1$ of the origine in $\cc_z$ 
$\pi^{-1}(z)$ is a disk;

\slii there is a neighborhood $W_2$ of  $[-1,1]$
on $\cc_z$ such that $\pi^{-1}(z)$ is an annulus for all $z\in W_2\setminus 
\bar W_1$;

\sliii for any $z\in [-1,1]\setminus \{0\}$ $\pi^{-1}(z)\cap f(H)$ is a 
circle $\{z\}\times \{|w|=1-\eps \}$ denoted by $\gamma_z$.

\smallskip Consider $(V,p)$ as a domain over $Y$ and let $(\tilde V,\tilde p)$
be its envelope of holomorphy. $\pi|_V $ extedns to a holomorphic function 
$\tilde \pi :\tilde V\to \Delta$ and it follows from the continuity principle
that $\tilde\pi^{-1}(z)$ contains a disk adjacent to the annulus $\pi^{-1}_z$
for all $z\in W_2$.

Furthermore, the $f_k$ holomorphically extends onto $\tilde V$. Therefore 
${\sf ind}_{\gamma_z}f_k\ge 0$, which contradicts Lemma 1.

\qed
\smallskip\noindent\bf Remark. \rm Using techniques from [Iv] one can 
show that $f(H)$ has no neighborhood biholomorphic to an open set in any 
holomorphically convex K\"ahler manifold.

\smallskip\noindent{\sl Counterexample to the Question 2.} 
Let now $K_j:|z-c_j|<r_j$ be a family
of mutually disjoint discs in $\Delta$ converging to $0$ and
$\Sigma_j$ be the intersection of $\bar K_j\times\Delta$ with
$\{r_j/3\le |z-c_j|\le 2r_j/3, w^2=({z-c_j\over
r_j})^{k_j}\lambda({|z-c_j|^2\over r_j^2})\}\cup \{|z-c_j|\ge r_j/3,
w=0\}$. Let $X$ be the domain $\Delta^2\setminus \left(\bigcup_j\Sigma_j
\cup\ \{ z\not\in \bigcup_jK_j,w=0\}\right)$
%(\cup K_j\times\Delta)\cup(\Delta\times A_{1/2})$
and $J$ be the complex structure (integrable!) in $X$ with the
basis of (1,0)-forms constituted by $dz$ and $dw+b\,d\bar z$ where
$b(z,w)=a_{k_j}({z-a_j\over r_j},w)$ for $z\in K_j,\ j=1,2,...$
and $b=0$ otherwise. Then $J\in C^\infty(X)$ if
$k_j\uparrow\infty$ sufficiently fast. $\pi :X\to \Delta$ denotes 
the natural projection which $X$ inherits as a domain in $\Delta\times
\Delta$.

Checking of the following Lemma is starightforward (due to Lemma 2) 
and is left to the reader.

\medskip
\begin{lem} Projection $\pi $ is holomorphic and therefore $(X,\pi)$ is a 
holomorphic fibration. Moreover:

\smallskip
\sli $X_z$ are disks with punctures; $X_0$ is a punctured disk; for 
$a\in \bigcup_j\{ |z-c_j|<r_j/3\}$ $X_a$ is a disk;

\smallskip
\slii there exists no foliated holomorphic map $(z,w):H_{\eps}\to (X,J)$
such that $|z(0)-c_j|<r_j/3$ for some $j$ and $z(1)=0$. 
\end{lem}

%\flushpar {\bf 2.} Another complex structure $J'$ in $M$ with (1,0)-basis
%$dz, dw+a\,d\bar w$ where $a=\eps_j\Lambda(|w|)$ in $K_j,\
%\Lambda\in C^\infty(\rr),\ 0\le\Lambda\le1,\ \Lambda(0)=1,\ \Lambda=0$
%for $t>1/2$ and $0<\eps_j<1$ are different numbers.
%Then $(M,J')$ admits holomorphic mappings from $(G,J_0)$ but a
%"natural hull" for such maps is non-Hausdorff $\tilde M$ "over" $M$
%obtained by glueing of countably many $\Delta^2$ by $\Delta\times A_{1/2}$.

\ifx\undefined\bysame
\newcommand{\bysame}{\leavevmode\hbox to3em{\hrulefill}\,}
\fi

\def\entry#1#2#3#4\par{\bibitem[#1]{#1}
{\textsc{#2 }}{\sl{#3} }#4\par\vskip2pt}
%{ref}{author}{title}.

\bigskip

Universite de Lille-I         \hfill Steklov Math. Institute

UFR de Mathematiques,         \hfill Russian Academi of Sci.

59655 Villeneuve d'Ascq,      \hfill Gubkin st. 8,

France                        \hfill 119991 Moscow, Russia

ivachkov@math.univ-lille1.fr  \hfill chirka@mi.ras.ru

\end{document}